\documentclass[12pt]{article}
\usepackage{algorithm,amsfonts,amsmath,amssymb,bbm,fullpage,graphicx,varwidth}
\usepackage{mathtools,psfrag,verbatim}
\usepackage[noend]{algpseudocode}
\usepackage[small,bf]{caption}

\newcommand{\reals}{{\mbox{\bf R}}}

\newcommand{\prox}{\mathbf{prox}}

\newcommand{\argmin}{\mathop{\rm argmin}}

\newcommand{\ie}{\emph{i.e.}}
\newcommand{\eg}{\emph{e.g.}}

\newcommand{\vnorm}[1]{\left\|#1\right\|}

\newcommand{\maxi}{\text{maximize}}
\newcommand{\mini}{\text{minimize}}

\title{Parameter Selection and Pre-Conditioning\\ for a Graph Form Solver}
\author{Christopher Fougner \and Stephen Boyd}

\begin{document}
\maketitle

\makeatletter
\def\BState{\State\hskip-\ALG@thistlm}
\makeatother

\algnewcommand\algorithmicinput{\textbf{Input:}}
\algnewcommand\Input{\item[\algorithmicinput]}

\begin{abstract}
In a recent paper, Parikh and Boyd describe a method for
solving a convex optimization problem,
where each iteration involves evaluating a proximal operator
and projection onto a subspace.
In this paper we address the critical practical issues
of how to select the proximal parameter in each iteration,
and how to scale the original problem variables, so as the 
achieve reliable practical performance.
The resulting method has been implemented as an open-source
software package called POGS (Proximal Graph Solver),
that targets multi-core and GPU-based systems,
and has been tested on a wide variety of practical problems.
Numerical results show that POGS
can solve very large problems (with, say, more than a billion coefficients
in the data), to modest accuracy in a few tens of seconds.
As just one example, a radiation treatment planning
problem with around 100 million coefficients in the data can be solved 
in a few seconds, as compared to around one hour with an interior-point method.
\end{abstract}

\section{Introduction} 

We consider the convex optimization problem
\begin{align}
	\begin{aligned}
	&\mini
	& & f(y) + g(x)  \\
	& \text{subject to} 
	& & y = A x,
	\end{aligned} \label{eq:gf_pri}
\end{align}
where $x\in \reals^n$ and $y \in \reals^{m}$ are the variables, and
the (extended-real-valued) functions $f: \reals^{m} \to \reals \cup \{\infty\}$ and $g: \reals^{n} \to \reals \cup \{\infty\}$ are convex, closed and proper. The matrix $A \in \reals^{m \times n}$, and the functions $f$ and $g$ are the problem data.
Infinite values of $f$ and $g$ allow us to encode convex 
constraints on $x$ and $y$, since any feasible point $(x,y)$ must satisfy
\[
x \in \{x \mid g(x) < \infty\}, \qquad
y \in \{y \mid f(y) < \infty\}.
\]
We will be interested in the case when $f$ and $g$ have simple proximal operators,
but for now we do not make this assumption.
The problem form~(\ref{eq:gf_pri}) 
is known as \emph{graph form} \cite{parikh2013block}, since the variable
$(x, y)$ is constrained to lie in the graph 
$\mathcal G = \{(x,y)  \in \reals^{n + m}~|~y = Ax\}$ of $A$.
We denote $p^\star$ as the optimal value of (\ref{eq:gf_pri}), which we assume
is finite.

The graph form includes a large range of convex problems, including
linear and quadratic programming, general conic programming \cite[\S 11.6]{BoV:04},
and many more 
specific applications such as logistic regression with various regularizers,
support vector machine fitting \cite{hastie2009elements},
portfolio optimization \cite[\S 4.4.1]{BoV:04} 
\cite{Glowinski1975} \cite{boyd2013performance},
and radiation treatment planning \cite{olafsson2006efficient}, to name just a few.

In \cite{parikh2013block}, Parikh and Boyd described an operator splitting 
method for solving the graph form problem~(\ref{eq:gf_pri}),
based on the alternating direction method of multipliers (ADMM)
\cite{boyd2011distributed}.
Each iteration of this method requires a projection 
(either exactly or approximately via an iterative method) onto the graph
$\mathcal G$, and evaluation of the proximal operators of $f$ and $g$.
Theoretical convergence was established in that paper, and basic
implementations demonstrated.
However it has been observed that practical convergence of the algorithm
depends very much on the choice of algorithm parameters 
(such as the proximal parameter $\rho$), and scaling of the variables
(\ie, pre-conditioning).

The purpose of this paper is to explore these issues, and to add
some critical variations on the algorithm that make it a relatively robust
general purpose solver, at least for modest accuracy levels.
The algorithm we propose, which is the same as the basic method 
described in \cite{parikh2013block}, with modified parameter selection,
diagonal pre-conditioning, and modified stopping criterion, has been 
implemented in an open-source software project called POGS
(for \textbf{P}r\textbf{o}ximal \textbf{G}raph \textbf{S}olver),
and tested on a wide variety of problems.
Our CUDA implementation reliably solves (to modest accuracy) 
problems $1000 \times$ larger than
those that can be handled by interior-point methods; and for those that
can be handled by interior-point methods, $100\times$ faster.
As a single example, a radiation treatment planning problem 
with more than 100 million
coefficients in $A$ can be solved in a few seconds; the same problem takes
around one hour to solve using an interior-point method.

\subsection{Outline} 
In \S\ref{sec:rel_work} we describe related work.
In \S\ref{sec:opt_dual} we derive the graph form dual problem,
and the primal-dual optimality conditions, which we use to motivate the
stopping criterion and to interpret the iterates of the algorithm.
In \S\ref{sec:alg} we describe the ADMM-based graph form algorithm, and
analyze the properties of its iterates, giving some results that
did not appear in \cite{parikh2013block}.
In \S\ref{sec:precond_param} we address the topic of
pre-conditioning, and suggest novel pre-conditioning and parameter selection
techniques.
In \S \ref{sec:impl} we describe our implementation POGS, and
in \S \ref{sec:numres} we report performance results on various problem families.

\subsection{Related work} \label{sec:rel_work}

Many generic methods can be used to solve the graph form 
problem~(\ref{eq:gf_pri}), including
projected gradient descent \cite{calamai1987projected}, projected 
subgradient methods
\cite[Chap. 5]{polyak1987introduction} \cite{shor1998nondifferentiable}, operator splitting methods
\cite{lions1979splitting} \cite{eckstein2008family}, interior-point methods
\cite[Chap. 19]{wrightnocedal1999numerical} 
\cite[Chap. 6]{ben2001lectures} and many more.
(Of course many of these methods can only be used when additional assumptions
are made on $f$ and $g$, \eg, differentiability or strong convexity.)
For example, if $f$ and $g$ are separable and
smooth (or have smooth barrier functions 
for their epigraphs), the problem (\ref{eq:gf_pri}) can be solved 
by an interior-point method, which in practice always takes no more 
than a few tens of iterations, with each iteration involving the solution of 
a system of linear equations that requires $O(\max\{m,n\}\min\{m,n\}^2)$
flops when $A$ is dense \cite[Chap. 11]{BoV:04}\cite[Chap. 19]{wrightnocedal1999numerical}.

We now turn to first-order methods for the graph form problem~(\ref{eq:gf_pri}).
In \cite{o2014primal} O'Connor and Vandenberghe
propose a primal-dual method for the graph form problem where $A$ is the sum of
two structured matrices. They contrast it with methods such as Spingarn's
method of partial inverses \cite{spingarn1985applications}, Douglas-Rachford
splitting \cite{douglas1956numerical}, and the Chambolle-Pock method
\cite{chambolle2011first}.

Davis and Yin \cite{davis2014convergence} analyze convergence rates for
different operator splitting methods, and in \cite{giselsson2015tight}
Giselsson proves the tightness of linear convergence for the operator
splitting problems considered \cite{giselsson2014metric}. Goldstein et al.
\cite{goldstein2014fast} derive Nesterov-type acceleration, and show
$O(1/k^2)$ convergence for problems where $f$ and $g$ are both
strongly convex.

Nishihara et al. \cite{nishihara2015general} introduce a parameter selection
framework for ADMM with over relaxation \cite{eckstein1992douglas}. The
framework is based on solving a fixed-size semidefinite program (SDP). They also
make the assumption that $f$ is strongly convex.  Ghadimi et al.
\cite{ghadimi2013optimal} derive optimal parameter choices for the case when
$f$ and $g$ are both quadratic.  In \cite{giselsson2014metric} Giselsson and
Boyd show how to choose metrics to optimize the convergence bound, and in
\cite{giselsson2014diagonal} Giselsson and Boyd suggest a diagonal
pre-conditioning scheme for graph form problems based on semidefinite
programming. This scheme is primarily relevant in small to medium scale
problems, or situations where many different graph form problems,
with the same matrix $A$, are to be solved.

It is clear from these papers (and indeed, a general rule) 
that the practical convergence of first-order mthods
depends heavily on algorithm parameter choices. 
All of these papers make additional assumptions about the objective,
which we do not.

GPUs are used extensively for training neural networks
\cite{ngiam2011optimization, ciresan2011flexible, krizhevsky2012imagenet,
coates2013deep}, and they are slowly gaining popularity in convex optimization
as well \cite{pock2011diagonal,chu2013primal,wang2014bregman}. 


\section{Optimality conditions and duality} \label{sec:opt_dual}

\subsection{Dual graph form problem}
The Lagrange dual function of (\ref{eq:gf_pri}) is given by
\[
	\inf_{x,y} f(y) + g(x) + \nu^T(Ax - y) = - f^*(\nu) -g^*(-A^T\nu)
\]
where $\nu\in \mathbf{R}^n$ is the dual variable associated with the equality 
constraint,
and $f^*$ and $g^*$ are the conjugate functions of $f$ and $g$ respectively
\cite[Chap. 4]{BoV:04}.
Introducing the variable $\mu = -A^T\nu$, we can write the dual
problem as
\begin{align}
	\begin{aligned}
	&\maxi
	& & -f^*(\nu) - g^*(\mu)  \\
	& \text{subject to} 
	& & \mu = -A^T \nu.
	\end{aligned} \label{eq:gf_dual}
\end{align}
The dual problem can be written as a graph form problem 
if we negate the objective and minimize rather than maximize.
The dual graph form problem~(\ref{eq:gf_dual})
is related to the primal graph form problem~(\ref{eq:gf_pri})
by switching the roles of the variables, replacing the objective function
terms with their conjugates, and replacing $A$ with $-A^T$.

The primal and dual objectives are $p(x,y) = f(y) + g(x)$ and $d(\mu,\nu) = -f^*(\nu) - g^*(\mu)$ respectively, giving us the duality gap
\begin{align}
	\eta = p(x,y) - d(\mu,\nu) = f(y)  + f^*(\nu) + g(x) + g^*(\mu). \label{eq:gap}
\end{align}
We have $\eta\geq 0$, for any primal and dual feasible tuple $(x, y, \mu, \nu)$.
The duality gap $\eta$ gives a bound on the suboptimality of $(x,y)$
(for the primal problem) and also $(\mu,\nu)$ for the dual problem:
\[
f(y)+g(x) \leq p^\star + \eta,  \qquad
-f^*(\nu)-g^*(\mu) \geq p^\star - \eta.
\]

\subsection{Optimality conditions}
The optimality conditions for~(\ref{eq:gf_pri}) 
are readily derived from the dual problem.
The tuple $(x, y, \mu, \nu)$ satisfies the following three conditions if and only
it is optimal.

\bigskip

\noindent \emph{Primal feasibility:}
\begin{align}
	y = Ax. \label{eq:pri_feas}
\end{align}
\emph{Dual feasibility:}
\begin{align}
	\mu = -A^T\nu.  \label{eq:dual_feas}
\end{align}
\emph{Zero gap:}
\begin{align}
	f(y)  + f^*(\nu) + g(x) + g^*(\mu) = 0. \label{eq:gap_opt}
\end{align}

If both (\ref{eq:pri_feas}) and (\ref{eq:dual_feas}) hold, 
then the zero gap condition (\ref{eq:gap_opt}) can be replaced by the Fenchel equalities
\begin{align}
	f(y) + f^*(\nu) = \nu^Ty, \quad g(x) + g^*(\mu) = \mu^Tx. \label{eq:fen_eq}
\end{align}
We refer to a tuple $(x,y,\mu,\nu)$ that satisfies~(\ref{eq:fen_eq})
as \emph{Fenchel feasible}.
To verify the statement, we add the two equations in~(\ref{eq:fen_eq}),
which yields
\begin{align*}
	f(y) + f^*(\nu) + g(x) + g^*(\mu) = y^T\nu + x^T\mu = (Ax)^T\nu -x^T A^T\nu = 0. 
\end{align*}
The Fenchel equalities (\ref{eq:fen_eq}) are is also equivalent to  
\begin{align}
	\nu \in \partial f(y), \quad \mu \in \partial g(x),  \label{eq:fen_sg} 
\end{align}
where $\partial$ denotes the subdifferential, which follows because
\[
	\nu \in \partial f(y) \Leftrightarrow \sup_z\left(z^T\nu - f(z) \right) =  \nu^Ty -  f(y) \Leftrightarrow f(y) + f^*(\nu) = \nu^Ty.
\]

In the sequel we will assume that strong duality holds, 
meaning that there exists a tuple $(x^\star, y^\star, \mu^\star, \nu^\star)$
which satisfies all three optimality conditions.

\section{Algorithm} \label{sec:alg}

\subsection{Graph projection splitting}
In \cite{parikh2013block} Parikh et al. apply ADMM \cite[\S 5]{boyd2011distributed}
to the problem of 
minimizing $f(y)+g(x)$, subject to the constraint $(x,y)\in \mathcal G$.
This yields the \emph{graph projection splitting} algorithm~\ref{alg:admm_pri}.

\begin{algorithm}
\caption{Graph projection splitting}
\label{alg:admm_pri}
\begin{algorithmic}[1]
\Input{$A, f, g$}
\State Initialize $(x^0, y^0, \tilde x^0, \tilde y^0) =0, ~k=0$
\Repeat
\State $(x^{k+1/2}, ~y^{k+1/2}) := \big(\mathbf{prox}_{g}(x^k-\tilde x^{k}),~ \mathbf{prox}_{f}(y^k-\tilde y^{k})\big)$
\State $(x^{k+1}, y^{k+1}) := \Pi(x^{k+1/2}+\tilde x^{k}, ~y^{k+1/2}+\tilde y^{k})$
\State $(\tilde x^{k+1}, \tilde y^{k+1}) := (\tilde x^{k} + x^{k+1/2} - x^{k+1}, ~\tilde y^{k} + y^{k+1/2} - y^{k+1})$
\State $k := k +1$
\Until {converged}
\end{algorithmic}
\end{algorithm}

The variable $k$ is the iteration counter, $x^{k+1}, x^{k+1/2} \in \mathbf{R}^{n}$ and 
$y^{k+1}, y^{k+1/2}, \in \mathbf{R}^{m}$ are primal variables,  $\tilde x^{k+1} \in \mathbf{R}^{n}$
and $\tilde y^{k+1} \in \mathbf{R}^{m}$ are scaled dual variables, $\Pi$ denotes the
(Euclidean) projection onto the graph $\mathcal G$,
\[
	\mathbf{prox}_{f}(v) = \argmin_{y}\Big(f(y) + (\rho/2) \vnorm{y - v}_2^2\Big)
\]
is the proximal operator of $f$ (and similarly for $g$),
and $\rho>0$ is the proximal parameter.
The projection $\Pi$ can be explicitly expressed as the linear operator
\begin{align}
	\Pi(c,d)= 
	K^{-1}\begin{bmatrix}
		c + A^Td \\ 0
	\end{bmatrix}, \qquad K = \begin{bmatrix}
		I & A^T \\
		A & -I
	\end{bmatrix}.
	\label{eq:proj}
\end{align}

Roughly speaking, in steps~3 and~5, the 
$x$ (and $\tilde x$) and $y$
(and $\tilde y$) variables do not mix;
the computations can be carried out in parallel.
The projection step~4 mixes the $x,\tilde x$ and $y,\tilde y$ variables.

General convergence theory for 
ADMM \cite[\S 3.2]{boyd2011distributed} guarantees that (with our assumption
on the existence of a solution)
\begin{align}
	(x^{k+1},y^{k+1}) - (x^{k+1/2}, y^{k+1/2}) \to 0, \quad f(y^k) +g(x^k) \to p^\star, \quad (\tilde x^k,\tilde y^k) \to (\tilde x^\star, \tilde y^\star),  \label{eq:conv_theo}
\end{align}
as $k \to \infty$. 

\subsection{Extensions}

We discuss three common extensions that 
can be used to speed up convergence in practice:
over-relaxation, approximate projection, and varying penalty.

\paragraph{Over-relaxation.} Replacing $x^{k+1/2}$ by $\alpha
x^{k+1/2}+(1-\alpha)x^{k}$ in the projection and dual update steps is known as
over-relaxation if $\alpha > 1$ or under-relaxation if $\alpha < 1$. 
The algorithm is guaranteed to converge 
\cite{eckstein1992douglas} for any $\alpha \in (0,2)$; 
it is observed in practice
\cite{o2013splitting} \cite{annergren2012admm} that
using an over-relaxation parameter in the range [1.5, 1.8]
can improve practical convergence.

\paragraph{Approximate projection.}
Instead of computing the projection $\Pi$ exactly one can use an approximation $\tilde \Pi$, with the only restriction that 
\[
	\textstyle{\sum}_{k = 0}^\infty\|\Pi(x^{k+1/2}, y^{k+1/2}) - \tilde \Pi(x^{k+1/2}, y^{k+1/2})\|_2 < \infty.
\]
This is known as approximate projection \cite{o2013splitting}. This extension is particularly useful if the approximate 
projection is computed using an indirect or iterative method.

\paragraph{Varying penalty.} Large values of $\rho$ tend to encourage primal
feasibility, while small values tend to encourage dual feasibility
\cite[\S3.4.1]{boyd2011distributed}. A common approach is to adjust or 
vary $\rho$ in
each iteration, so that the primal and dual residuals are (roughly) balanced in
magnitude.
When doing
so, it is important to re-scale $(\tilde x^{k+1}, \tilde y^{k+1})$ by a factor
$\rho^{k}/\rho^{k+1}$.

\subsection{Feasible iterates} \label{sec:feas_vars}
In each iteration, algorithm~\ref{alg:admm_pri} produces sets of points that are either primal, dual, or Fenchel feasible. Define
\begin{align*}
	\mu^{k} = -\rho\tilde x^k, \quad \nu^{k} = -\rho\tilde y^k, \quad  \mu^{k+1/2} = -\rho(x^{k+1/2} - x^k + \tilde x^k), \quad \nu^{k+1/2} = -\rho (y^{k+1/2} - y^k + \tilde y^k). 
\end{align*}
The following statements hold.
\begin{enumerate}
\item The pair $(x^{k+1}, y^{k+1})$  is primal feasible, 
since it is the projection onto the graph $\mathcal G$.
\item  \label{enum:fen} The pair $(\mu^{k+1}, \nu^{k+1})$ is dual feasible, as long as $(\mu^{0}, \nu^{0})$ is dual feasible and $(x^0, y^0)$ is primal feasible. Dual feasibility implies $\mu^{k+1} +A^T \nu^{k+1}=0$, which we show using the update equations in algorithm~\ref{alg:admm_pri}:
\begin{align*}
	\mu^{k+1} +A^T \nu^{k+1} &=  -\rho(\tilde x^{k} + x^{k+1/2} - x^{k+1}  +A^T ( \tilde y^{k} + y^{k+1/2} - y^{k+1} )) \\
		&=  -\rho(\tilde x^{k} + A^T\tilde y^k + x^{k+1/2} +A^Ty^{k+1/2} - (I+A^TA)x^{k+1}),
\end{align*}
where we substituted $y^{k+1} = Ax^{k+1}$. From the projection operator in (\ref{eq:proj}) it follows that  $(I+A^TA)x^{k+1} = x^{k+1/2} + A^Ty^{k+1/2}$, therefore
\[
	\mu^{k+1} +A^T \nu^{k+1}  = -\rho(\tilde x^{k} + A^T\tilde y^k) =  \mu^{k} +A^T \nu^{k} = \mu^0 + A^T\nu^0,
\]
where the last equality follows from an inductive argument. Since we made the assumption that $ (\mu^{0},\nu^{0})$ is dual feasible, we can conclude that $(\mu^{k+1}, \nu^{k+1})$ is also dual feasible.

\item The tuple $(x^{k+1/2}, y^{k+1/2}, \mu^{k+1/2}, \nu^{k+1/2})$ is Fenchel feasible. From the definition of the proximal operator,
\begin{align*}
	x^{k+1/2} =  \argmin_{x}\Big(g(x) + (\rho/2) \vnorm{x - x^k + \tilde x^k}_2^2\Big) &\Leftrightarrow  0 \in \partial g(x^{k+1/2}) + \rho(x^{k+1/2} - x^k + \tilde x^k) \\
	&\Leftrightarrow \mu^{k+1/2} \in \partial g(x^{k+1/2}).
\end{align*}
By the same argument $ \nu^{k+1/2} \in \partial f(y^{k+1/2})$.

\end{enumerate}

Applying the results in (\ref{eq:conv_theo}) to the dual variables, we find $\nu^{k+1/2} \to \nu^\star $ and $\mu^{k+1/2} \to \mu^\star$, from which we conclude that $(x^{k+1/2}, y^{k+1/2}, \mu^{k+1/2}, \nu^{k+1/2})$ is primal and dual feasible in the limit.

\subsection{Stopping criteria}

In \S\ref{sec:feas_vars} we noted that either (\ref{eq:pri_feas}, \ref{eq:dual_feas}, \ref{eq:gap_opt}) or (\ref{eq:pri_feas}, \ref{eq:dual_feas}, \ref{eq:fen_eq}) are sufficient for optimality. We present two different stopping criteria based on these conditions.

\paragraph{Residual based stopping.} \label{sec:res_stop}
The tuple $(x^{k+1/2}, y^{k+1/2}, \mu^{k+1/2}, \nu^{k+1/2})$ is Fenchel feasible in each iteration, but only primal and dual feasible in the limit. Accordingly, we propose the residual based stopping criterion
\begin{align}
	\|Ax^{k+1/2} - y^{k+1/2}\|_2 \leq \epsilon^{\text{pri}}, \quad \|A^T\nu^{k+1/2} + \mu^{k+1/2}\|_2 \leq \epsilon^{\text{dual}}, \label{eq:res_stop}
\end{align}
where the $\epsilon^\text{pri}$ and $\epsilon^\text{dua}$ are positive tolerances. These should be chosen as a mixture of absolute and relative tolerances, such as
\[
	\epsilon^\text{pri} = \epsilon^{\text{abs}} + \epsilon^{\text{rel}} \|y^{k+1/2}\|_2, \quad \epsilon^\text{dual} = \epsilon^{\text{abs}} + \epsilon^{\text{rel}} \|\mu^{k+1/2}\|_2.
\]
Reasonable values for $\epsilon^{\text{abs}}$ and $\epsilon^{\text{rel}}$ are in the range $[10^{-4}, 10^{-2}]$.

\paragraph{Gap based stopping.}
The tuple $(x^k, y^k, \mu^k, \nu^k)$ is primal and dual feasible, but only Fenchel feasible in the limit. We propose the gap based stopping criteria
\[
	\eta^k = 
f(y^k) + g(x^k) + f^*(\nu^k) + g^*(\mu^k) \leq \epsilon^\text{gap},
\]
where $\epsilon^\text{gap}$ should be chosen relative to the current objective value, \ie,
\[
	\epsilon^\text{gap} = \epsilon^{\text{abs}} + \epsilon^{\text{rel}}|f(y^k) + g(x^k)|.
\]
Here too, reasonable values for $\epsilon^{\text{abs}}$ and 
$\epsilon^{\text{rel}}$ are in the range $[10^{-4}, 10^{-2}]$. 

Although the gap based stopping criteria is very informative, since it directly 
bounds the suboptimality of the current iterate,
it suffers from the drwaback that 
$f, g, f^*$ and $g^*$ must all have full domain, since otherwise
the gap $\eta^k$ can be infinite.  
Indeed, the gap $\eta^k$ is almost always infinite
when $f$ or $g$ represent constraints.

\subsection{Implementation} \label{sec:impl_cons}
\paragraph{Projection.}
There are different ways to evaluate the projection operator $\Pi$,
depending on the structure and size of $A$. 

One simple method that can be used if $A$ is sparse and not too large is a
direct sparse factorization. The matrix $K$ is quasi-definite, and therefore
the $LDL^T$ decomposition  is well defined \cite{vanderbei1995symmetric}. Since
$K$ does not change from iteration to iteration, the factors $L$ and $D$ (and
the permutation or elimination ordering) can be computed in the first iteration
(\eg, using CHOLMOD \cite{chen2008algorithm}) and re-used in subsequent
iterations. This is known as \emph{factorization caching} 
\cite[\S 4.2.3]{boyd2011distributed} \cite[\S A.1]{parikh2013block}.
With factorization caching, we get a (potentially) large speedup in
iterations, after the first one.

If $A$ is dense, and $\min(m,n)$ is not too large,
then block elimination \cite[Appendix C]{BoV:04} can be applied to $K$ 
\cite[Appendix A]{parikh2013block}, yielding the reduced update
\begin{align*}
	x^{k+1} &:= (A^TA + I)^{-1}(c + A^Td) \\
	y^{k+1} &:= Ax^{k+1}
\end{align*}
if $m \geq n$, or  
\begin{align*}
	y^{k+1} &:= d + (AA^T + I)^{-1}(Ac -d) \\
	x^{k+1} &:= c - A^T(d-y^{k+1})
\end{align*}
if $m < n$. 
Both formulations involve forming and solving a system of equations 
in $\mathbf{R}^{\text{min}(m,n) \times \text{min}(m,n)}$. 
Since the matrix is symmetric positive definite, 
we can use the Cholesky decomposition.
Forming the coefficient matrix $A^TA+I$ or $AA^T+I$ dominates the computatation.
Here too we can take advantage of factorization caching.

The regular structure of dense matrices allows us to analyze the computational
complexity of each step. We define $q =\min(m,n)$ and $p = \max(m,n)$.  The
first iteration involves the factorization and the solve step; subsequent
iterations only require the solve step. The computational cost of the
factorization is the combined cost of computing $A^TA$ (or $AA^T$, whichever is
smaller), at a cost of $pq^2$ flops, in addition to the Cholesky decomposition,
at a cost of $(1/3)q^3$ flops. The solve step consists of two matrix-vector
multiplications at a cost of $4pq$ flops and solving a triangular system of
equations at a cost of $q^2$ flops. The total cost of the first iteration is
$O(pq^2)$ flops, while each subsequent iteration only costs $O(pq)$ flops,
showing that we obtain a savings by a factor of $q$ flops, after the first
iteration, by using factorization caching.

For very large problems direct methods are no longer practical, 
at which point indirect (iterative) methods can be used.
Fortunately, as the primal and dual variables converge, we are
guaranteed that $(x^{k+1/2}, y^{k+1/2}) \to (x^{k+1}, y^{k+1})$, meaning that
we will have a good initial guess we can use to 
initialize the iterative method to 
(approximately) evaluate the projection.
One can either apply CGLS (conjugate gradient least-squares) \cite{hestenes1952methods} or
LSQR \cite{paige1982lsqr} to the reduced update or apply MINRES (minimum residual)
\cite{paige1975solution} to $K$ directly. It can be shown the latter requires twice the number of iterations as compared to the former, and is therefore not recommended.

\paragraph{Proximal operators.}
Since the $x,\tilde x$ and $y,\tilde y$ components are decoupled in the proximal step and dual
variable update step, both of these can be done separately, and in parallel for
$x$ and $y$. If either $f$ or $g$ is separable, then the proximal step can be
parallelized further. The monograph \cite{parikh2013proximal} details how
proximal operators can be computed efficiently for a wide range of functions.
Typically the cost of computing the proximal operator will be negligible
compared to the cost of the projection. In particular, if $f$ and $g$ are
separable, then the cost will be $O(m + n)$, and completely parallelizable.

\section{Pre-conditioning and parameter selection} \label{sec:precond_param}

The practical convergence of the algorithm (\ie, the number of iterations required 
before it terminates) can depend greatly on the choice of the proximal 
parameter $\rho$, and the scaling of the variables.
In this section we analyze these, and suggest a method for choosing
$\rho$ and for scaling the variables that (empirically) 
speeds up practical convergence.

\subsection{Pre-conditioning} \label{sec:precond}

Consider scaling the variables $x$ and $y$ in (\ref{eq:gf_pri}), by $E^{-1}$ and $D$ respectively, where $D \in \mathbf{R}^{m\times m}$ and $E \in \mathbf{R}^{n \times n}$ are non-singular matrices. 
We define the scaled variables
\[
	\hat y = Dy, \quad \hat x = E^{-1}x,
\]
which transforms (\ref{eq:gf_pri}) into
\begin{align}
	\begin{aligned}
	&\mini
	& & f(D^{-1}\hat y) + g(E \hat x)  \\
	& \text{subject to} 
	& & \hat y= DAE \hat x. \label{eq:gf_trans}
	\end{aligned}
\end{align}
This is also a graph form problem, and for notational convenience, we define 
\[
	\quad \hat A = DAE, \quad \hat f(\hat y) = f(D^{-1}\hat y), \quad \hat g(\hat x) = g(E \hat x),
\]
so that the problem can be written as
\begin{align*}
	\begin{aligned}
	&\mini
	& & \hat f(\hat y) + \hat g(\hat x)  \\
	& \text{subject to} 
	& & \hat y = \hat A\hat x.
	\end{aligned}
\end{align*}
We refer to this problem as the pre-conditioned version of (\ref{eq:gf_pri}).
Our goal is to choose $D$ and $E$ so that (a) the algorithm applied to the 
pre-conditioned problem converges in fewer steps in practice, and (b) the 
additional computational cost due to the pre-conditioning is minimal.

Graph projection splitting applied to the pre-conditioned problem (\ref{eq:gf_trans}) can be interpreted in terms of the original iterates. The proximal step iterates are redefined as
\begin{align*}
	x^{k+1/2} &= \argmin_{x} \left( g(x) + (\rho/2)\|x - x^k + \tilde x^k\|_{(EE^T)^{-1}}^2 \right) \\
	y^{k+1/2} &= \argmin_{y} \left(f(y) + (\rho/2)\|y - y^k + \tilde y^k\|_{(D^TD)}^2 \right),
\end{align*}
and the projected iterates are the result of the weighted projection
\begin{align*}
	\begin{aligned}
	&\mini
	& & (1/2)\|x - x^{k+1/2}\|_{(EE^T)^{-1}}^2 + (1/2)\|y - y^{k+1/2}\|_{(D^TD)}^2  \\
	& \text{subject to} 
	& & y = A x,
	\end{aligned}
\end{align*}
where $\|x\|_P = \sqrt{x^TPx}$ for a symmetric positive-definite matrix $P$. This projection can be expressed as
\[
	\Pi(c,d)= 
	\hat K^{-1}\begin{bmatrix}
		(EE^T)^{-1}c + A^TD^TDd \\ 0
	\end{bmatrix}, \qquad \hat K = \begin{bmatrix}
		(EE^T)^{-1} & A^TD^TD \\
		D^TDA & -D^TD
	\end{bmatrix}.
\]

Notice that graph projection splitting is invariant to orthogonal
transformations of the variables $x$ and $y$, since the pre-conditioners only
appear in terms of $D^TD$ and $EE^T$. In particular, if we let $D = U^T$ and $E
= V$, where $A=U\Sigma V^T$, then the pre-conditioned constraint matrix $\hat A
= DAE  = \Sigma$ is diagonal.  We conclude that
any graph form problem can be pre-conditioned to one
with a diagonal non-negative constraint matrix $\Sigma$. For
analysis purposes, we are therefore free to assume that $A$ is diagonal. We also note that for orthogonal pre-conditioners, there exists an analytical relationship between the original proximal operator and the pre-conditioned proximal operator. With $\phi(x) =
\varphi(Qx)$, where $Q$ is any orthogonal matrix ($Q^TQ = QQ^T= I$), we have
\[
\mathbf{prox}_{\phi}(v) = Q^T\mathbf{prox}_{\varphi}(Qv).
\] 
While the proximal operator of $\phi$ is readily computed, orthogonal pre-conditioners destroy separability of the objective. As a result, we can not easily combine them with other pre-conditioners.

Multiplying $D$ by a scalar $\alpha$  and dividing $E$ by the same scalar has the effect of scaling $\rho$ by a factor of $\alpha^2$. It however has no effect on the projection step, showing that $\rho$ can be thought of as the relative scaling of $D$ and $E$.

In the case where $f$ and $g$ are separable and both $D$ and $E$ are diagonal, the proximal step takes the simplified form
\begin{align*}
	x_j^{k+1/2} &= \argmin_{x_j} \left(g_j(x_j) + (\rho^E_{j}/2)(x_j - x_j^k + \tilde x_j^k)^2 \right) ~&& j = 1,\ldots,n\\
	y_i^{k+1/2} &= \argmin_{y_i} \left(f_i(y_i) + (\rho^D_{i}/2)(y_i - y_i^k + \tilde y_i^k)^2 \right) ~&& i = 1,\ldots,m,
\end{align*}
where $\rho^E_{j} = \rho/E_{jj}^2$ and $\rho^D_{i} = \rho D_{ii}^2$. Since only $\rho$ is modified, any routine capable of computing $\mathbf{prox}_f$ and $\mathbf{prox}_g$ can also be used to compute the pre-conditioned proximal update.

\subsubsection{Effect of pre-conditioning on projection}

For the purpose of analysis, we will assume that $A = \Sigma$, where $\Sigma$ is a non-negative diagonal matrix. The projection operator simplifies to
\[
	\Pi(c,d) = \begin{bmatrix}
		(I + \Sigma^T\Sigma)^{-1} & (I + \Sigma^T\Sigma)^{-1}\Sigma^T \\
		(I + \Sigma\Sigma^T)^{-1}\Sigma &(I + \Sigma\Sigma^T)^{-1}\Sigma\Sigma^T
	\end{bmatrix}\begin{bmatrix} c \\ d\end{bmatrix},
\]
which means the projection step can be written explicitly as
\begin{align*}
	x^{k+1}_i &= \frac{1}{1+\sigma_i^2}(x_i^{k+1/2} + \tilde x_i^{k} + \sigma_i(y_i^{k+1/2} +  \tilde y_i^{k})) && 1 \leq i \leq \min(m,n) \\
	x^{k+1}_i &= x_i^{k+1/2} + \tilde x_i^{k}&& \min(m,n) < i \leq n \\
	y^{k+1}_i &= \frac{\sigma_i}{1+\sigma_i^2}(x_i^{k+1/2} + \tilde x_i^{k} + \sigma_i (y_i^{k+1/2} + \tilde y_i^{k})) && 1 \leq i \leq \min(m,n) \\
	y^{k+1}_i &= 0 && \min(m,n) < i \leq m,
\end{align*}
where $\sigma_i$ is the $i$th diagonal entry of $\Sigma$ and subscripted indices of $x$ and $y$ denote the $i$th entry of the respective vector. Notice that the projected variables $x_i^{k+1}$ and $y_i^{k+1}$ are equally dependent on $(x_i^{k+1/2} + \tilde x_i^k)$ and $\sigma_i(y_i^{k+1/2} +  \tilde y_i^{k})$. If $\sigma_i$ is either significantly smaller or larger than 1, then the terms $x_i^{k+1}$ and $y_i^{k+1}$ will be dominated by either $(x_i^{k+1/2} + \tilde x_i^k)$ or $(y_i^{k+1/2} +  \tilde y_i^{k})$. However if $\sigma_i = 1$, then the projection step exactly 
averages the two quantities
\begin{align*}
	x_i^{k+1} = y_i^{k+1} &= \frac{1}{2}(x_i^{k+1/2} + \tilde x_i^{k} + y_i^{k+1/2} + \tilde y_i^k)  && 1 \leq i \leq \min(m,n).
\end{align*}
As we pointed out in \S\ref{sec:alg}, the projection step mixes 
the variables $x$ and $y$.  For this to approximately reduce to averaging,
we need $\sigma_i \approx 1$.

\subsubsection{Choosing $D$ and $E$}

The analysis suggests that the algorithm will be fast when 
the singular values of $DAE$ are all near one, \ie,
\begin{align}
	\mathbf{cond}\big(DAE\big) \approx 1, \quad \|DAE\|_2 \approx 1. \label{eq:precond_obj}
\end{align}
(This claim is also supported in \cite{giselsson2014preconditioning},
and is consistent with our computational experience.)
Pre-conditioners that exactly satisfy these conditions can be found using
the singular value decomposition of $A$.
They will however be of little use, since such
pre-conditioners generally destroy our ability to evaluate the proximal operators
of $\hat f$ and $\hat g$ efficiently.

So we seek choices of $D$ and $E$ for which~(\ref{eq:precond_obj}) holds
(very) approximately, and for which the proximal operators of
$\hat f$ and $\hat g$ can still be efficiently computed.
We now specialize to the special case when $f$ and $g$ are 
separable.  In this case, diagonal $D$ and $E$ are 
candidates for which the proximal operators are still 
easily computed.  (The same ideas apply to block separable $f$ and $g$,
where we impose the further constraint that the diagonal entries within
a block are the same.)
So we now limit ourselves to the case of diagonal pre-conditioners.

Diagonal matrices that minimize the condition number of $DAE$, 
and therefore approximately satisfy the first condition in (\ref{eq:precond_obj}),
can be found exactly, using semidefinite programming \cite[\S 3.1]{boyd1994linear}.
But this computation is quite involved, and may not be worth the computational
effort since the conditions~(\ref{eq:precond_obj}) are just a heuristic 
for faster convergence.
(For control problems, where the problem is solved many times with the same
matrix $A$, this approach makes sense; see
\cite{giselsson2014diagonal}.)

A heuristic that tends to minimize the condition number is to equilibrate the
matrix, \ie, choose $D$ and $E$ so that the rows all have the same $p$-norm,
and the columns all have the same $p$-norm.  (Such a matrix is said to be
equilibrated.)
This corresponds to finding $D$ and $E$
so that 
\[
|DAE|^p\mathbf{1} = \alpha \mathbf{1}, \qquad \mathbf{1}^T|DAE|^p = \beta
\mathbf{1}^T,
\]
where $\alpha, \beta > 0$. Here the notation
$|\cdot|^p$ should be understood in the elementwise sense. Various authors
\cite{o2013splitting}, \cite{chu2013primal}, \cite{bradley2010algorithms}
suggest that equilibration can decrease the number of iterations needed for
operator splitting and other first order methods.
One issue that we need to address is that not every matrix can be equilibrated.
Given that equilibration is only a heuristic for achieving
$\sigma_i(DAE)\approx 1$, which is in turn a heuristic for
fast convergence of the algorithm, partial equilibration should
serve the same purpose just as well.

Sinkhorn and Knopp \cite{sinkhorn1967concerning} suggest a method
for matrix equilibration for $p<\infty$, and $A$ is square and has full
support. In the case $p = \infty$, the Ruiz algorithm
\cite{ruiz2001scaling} can be used.   Both of these methods fail (as 
they must) when the matrix $A$ cannot be equilibrated.
We give below a simple modification of the
Sinkhorn-Knopp algorithm, 
modified to handle the case when $A$ is non-square, or cannot be 
equilibrated.

Choosing pre-conditioners that satisfy $\|DAE\|_2=1$ can be achieved by scaling
$D$ and $E$ by $\sigma_{\max}(DAE)^{-q}$ and $\sigma_{\max}(DAE)^{q-1}$
respectively for $q \in \mathbf{R}$. The quantity $\sigma_{\max}(DAE)$
can be approximated using power iteration, but we have found
it is unnecessary to exactly enforce $\|DAE\|_2=1$.
A more computationally efficient alternative
is to replace $\sigma_{\max}(DAE)$ by $\|DAE\|_F/\sqrt{\min(m,n)}$. This
quantity coincides with $\sigma_{\max}(DAE)$ when $\mathbf{cond}(DAE) = 1$. If
$DAE$ is equilibrated and $p=2$, this scaling corresponds to $(DAE)^T(DAE)$ (or
$(DAE)(DAE)^T$ when $m < n$) having unit diagonal.

\subsection{Regularized equilibration}

In this section we present a self-contained derivation of 
our matrix-equilibration
method. It is similar to the Sinkhorn-Knopp algorithm, but also works when the
matrix is non-square or cannot be exactly equilibrated.

Consider the convex optimization problem with variables $u$ and $v$,
\begin{align}
	\begin{aligned}
	&\mini
	&  & \sum_{i=1}^m\sum_{j=1}^n|A_{ij}|^pe^{u_i+v_j} - n \mathbf{1}^Tu - m \mathbf{1}^Tv +\gamma \left[(1/m)\sum_{i=1}^me^{u_i} + (1/n)\sum_{j=1}^ne^{v_j}\right], \label{eq:sk_obj}
	\end{aligned}
\end{align}
where $\gamma \geq 0$ is a regularization parameter. The objective is bounded below for any $\gamma > 0$.
The optimality conditions are
\begin{align*}
	\sum_{j=1}^n|A_{ij}|^pe^{u_i+v_j} - n + (1/m)\gamma e^{u_i}= 0,\quad i = 1,\ldots,m\\
	\sum_{i=1}^m|A_{ij}|^pe^{u_i+v_j} - m + (1/n)\gamma e^{v_j} = 0,\quad j = 1,\ldots,n.
\end{align*}
By defining $D_{ii} = e^{u_i/p}$ and $E_{jj}^p = e^{v_j/p}$, these conditions are equivalent to
\[
	|DAE|^p \mathbf{1} + (1/m)\gamma D\mathbf{1} =n \mathbf{1}, \quad \mathbf{1}^T |DAE|^p  +  (1/n)\gamma \mathbf{1}^TE=m \mathbf{1}^T.
\]
When $\gamma = 0$, these are the conditions for a matrix to be equilibrated. The objective may not be bounded when $\gamma = 0$, which exactly corresponds to the case when the matrix cannot be equilibrated. As $\gamma \to \infty$, both $D$ and $E$ converge to the scaled identity matrix $(mn/\gamma)I$, showing that $\gamma$ can be thought of as a regularizer on the elements of $D$ and $E$. If $D$ and $E$ are optimal, then the two equalities
\[
	\mathbf{1}^T|DAE|^p\mathbf{1} +  (1/m)\gamma \mathbf{1}^TD\mathbf{1}  =  mn,
\qquad \mathbf{1}^T|DAE|^p\mathbf{1} +  (1/n)\gamma \mathbf{1}^TE\mathbf{1} = mn
\]
must hold. Subtracting the one from the other, and dividing by $\gamma$, we find the relationship
\[
	(1/m) \mathbf{1}^TD\mathbf{1} = (1/n) \mathbf{1}^TE\mathbf{1},
\]
implying that the average entry in $D$ and $E$ is the same.

There are various ways to solve the optimization problem (\ref{eq:sk_obj}), one of which is to apply coordinate descent. Minimizing the objective in 
(\ref{eq:sk_obj}) with respect to $u_i$ yields
\[
	\sum_{j = 1}^ne^{u_i^{k}+v^{k-1}_j}|A_{ij}|^p + (\gamma/m) e^{u_i^{k}}= n \Leftrightarrow e^{u_i^{k}} = \frac{n}{\sum_{j = 1}^ne^{v^{k-1}_j}|A_{ij}|^p + (\gamma/m)}
\]
and equivalently for $v_j$
\[
	e^{v_i^{k}} = \frac{m}{\sum_{i = 1}^ne^{u^{k-1}_i}|A_{ij}|^p + (\gamma/n)}.
\]
Since the minimization with respect to $u_i^k$ is independent of $u_{i-1}^{k}$, the update can be done in parallel for each element of $u$, and similarly for $v$. 
Repeated minimization over $u$ and $v$ will eventually yield values that 
satisfy the optimality conditions.  Algorithm~\ref{alg:sk1} summarizes the equilibration routine.
\begin{algorithm}
\caption{Regularized Sinkhorn-Knopp}
\label{alg:sk1}
\begin{algorithmic}[1]
\Input{$A, \epsilon>0$, $\gamma>0$}
\State Initialize $e^0 := \mathbf{1}, ~k := 0$
\Repeat
\State $k := k +1$
\State $d^{k} := n ~ \mathbf{diag}(|A|^pe^{k-1} + (\gamma/m)\mathbf{1})^{-1} \mathbf{1}$
\State $e^{k} := m ~ \mathbf{diag}(|A^T|^pd^{k} + (\gamma/n)\mathbf{1})^{-1} \mathbf{1}$
\Until {$\|e^k - e^{k-1}\|_2 \leq \epsilon$ and $\|d^k - d^{k-1}\|_2 \leq \epsilon$}
\State \Return $D := \mathbf{diag}(d^{k})^{1/p}$, $E := \mathbf{diag}(e^{k})^{1/p}$
\end{algorithmic}
\end{algorithm}

\subsection{Adaptive penalty update}

The projection operator $\Pi$ does not depend on the choice of $\rho$, so we are
free to update $\rho$ in each iteration, at no extra cost. While the
convergence theory only holds for fixed $\rho$, it still applies if one assumes
that $\rho$ becomes fixed after a finite number of iterations
\cite{boyd2011distributed}. 

As a rule, increasing $\rho$ will decrease the
primal residual, while decreasing $\rho$ will decrease the dual residual. 
The authors in \cite{he2000alternating},\cite{boyd2011distributed} 
suggest adapting $\rho$ to balance the
primal and dual residuals.  We have found that substantially better 
practical convergence can be obtained using a variation on this idea.
Rather than balancing the primal and dual residuals, we allow
either the primal or dual residual to approximately converge and only
then start adjusting $\rho$. 
Based on this observation, we propose the following adaptive
update scheme.

\begin{algorithm}
\caption{Adaptive $\rho$ update}
\label{alg:adaptive_rho}
\begin{algorithmic}[1]
\Input{$\delta > 1, \, \tau \in (0,1], $}
\State Initialize $l:= 0, \,u := 0$
\Repeat 
\State Apply graph projection splitting
\If {$\|A^T\nu^{k+1/2}+\mu^{k+1/2}\|_2 < \epsilon^{\text{dual}}$ and $\tau k > l $}
   \State $\rho^{k+1} := \delta \rho^k$
   \State $u := k$
\ElsIf {$\|Ax^{k+1/2} - y^{k+1/2}\|_2 < \epsilon^{\text{pri}}$  and $\tau k > u $}
   \State $\rho^{k+1} := (1/\delta)\rho^k$
   \State $l := k$
\EndIf
\Until {$\|A^T\nu^{k+1/2}+\mu^{k+1/2}\|_2 < \epsilon^{\text{dual}}$  and  $\|Ax^{k+1/2} - y^{k+1/2}\|_2 < \epsilon^{\text{pri}}$  }
\end{algorithmic}
\end{algorithm}
Once either the primal or dual residual converges, the algorithm begins to steer $\rho$ in a direction so that the other residual also converges. 
By making small adjustments to $\rho$, we will tend to remain approximately 
primal (or dual) feasible once primal (dual) feasibility has been attained.
Additionally by requiring a certain number of iterations between an increase in $\rho$ and a decrease (and vice versa),  we enforce that changes to $\rho$ do not flip-flop between one direction and the other. The parameter $\tau$ determines the relative number of iterations between changes in direction. 

\section{Implementation} \label{sec:impl}

Proximal Graph Solver (POGS) is an open-source (BSD-3 license) implementation of graph projection splitting, written in C++. It supports both GPU and CPU platforms and includes wrappers for C, MATLAB, and R. POGS handles all combinations of sparse/dense matrices, single/double precision arithmetic, and direct/indirect solvers, with the exception (for now)
of sparse indirect solvers. The only dependency is a tuned BLAS library on the respective platform (\eg, cuBLAS or the Apple Accelerate Framework). The source code is available at
\begin{center}
\begin{varwidth}{\linewidth}
\small \begin{verbatim}
https://github.com/cvxgrp/pogs
\end{verbatim} \normalsize
\end{varwidth}
\end{center}

In lieu of having the user specify the proximal operators of $f$ and $g$, POGS contains a library of proximal operators for a variety of different functions. It is currently assumed that the objective is separable, in the form
\[
	f(y) + g(x) = \sum_{i = 1}^mf_i(y_i) + \sum_{j = 1}^ng_j(x_j),
\]
where $f_i, g_j : \mathbf{R} \to \mathbf{R} \cup \{\infty\}$. The library contains a set of base functions, and by applying various transformations, the range of functions can been greatly extended. In particular we use the parametric representation 
\[
	f_i(y_i) = c_i h_i(a_iy_i - b_i) + d_iy_i + (1/2)e_iy_i^2,
\]
where $a_i,b_i,d_i \in \mathbf{R}$, $c_i, e_i \in \mathbf{R}_+$, and $h_i : \mathbf{R} \to \mathbf{R} \cup \{\infty\}$. The same representation is also used for $g_j$. It is straightforward to express the proximal operators of $f_i$ in terms of the proximal operator of $h_i$ using the formula
\begin{align*}
	{\prox}_{f}(v) = \frac{1}{a}\bigg({\prox}_{h, (e+\rho)/(c a^2)}\Big(a \left(v\rho -d\right)/(e+\rho) - b\Big) + b\bigg),
\end{align*}
where for notational simplicity we have dropped the index $i$ in the constants 
and functions.
It is possible for a user to add their own proximal operator function, if
it is not in the current library.
We note that the separability assumption on $f$ and $g$ is a simplification,
rather than a limitation of the algorithm.
It allows us to apply the proximal operator in
parallel using either CUDA or OpenMP (depending on the platform). 

The constraint matrix is equilibrated using algorithm~\ref{alg:sk1}, with a
choice of $p=2$ and $\gamma = (m + n)\sqrt{\epsilon^{\text{cmp}}}$, where
$\epsilon^{\text{cmp}}$ is machine epsilon. Both $D$ and $E$ are rescaled
evenly, so that they satisfy $\|DAE\|_F/\sqrt{\min(m,n)} = 1$.  The projection
$\Pi$ is computed as outlined in \S\ref{sec:impl_cons}. We work with the
reduced update equations in all versions of POGS. In the indirect case, we chose
to use CGLS.
The parameter $\rho$ is updated according to algorithm~\ref{alg:adaptive_rho}.
Empirically, we found that $(\delta, \,\tau) = (1.05, \,0.8)$ works well. We
also use over-relaxation with $\alpha = 1.7$. 

POGS supports warm starting, whereby an initial guess for $x^0$ and/or $\nu^0$ may be supplied by the user.
If only $x^0$ is provided, then $\nu^0$ will be estimated, and vice-versa. 
The warm-start feature allows any cached matrices to be used 
to solve additional problems with the same matrix $A$.

POGS returns the tuple ($x^{k+1/2}, y^{k+1/2}, \mu^{k+1/2}, \nu^{k+1/2}$), since
it has finite primal and dual objectives. The primal and dual residuals will be
non-zero and are determined by the specified tolerances.

Future plans for POGS include extension to block-separable $f$ and $g$ (including 
general cone solvers), additional
wrappers for Julia and Python, support for a sparse direct solver, and
a multi-GPU extension. 

\section{Numerical results} \label{sec:numres}

To highlight the robustness and general purpose nature of POGS, we tested it on
9 different problem classes using random data, as well as a
radiation treatment planning problem using real-world data.

All experiments were performed in single precision arithmetic on a machine
equipped with an Intel Core i7-870, 16GB of RAM, and a Tesla K40 GPU. Timing
results include the data copy from CPU to GPU. 

We compare POGS to SDPT3 \cite{toh1999sdpt3}, an open-source solver
that handles linear, second-order, and positive semidefinite cone programs.
Since SDPT3 uses an interior-point algorithm, the solution returned will be
of high precision, allowing us to verify the accuracy of the solution computed
by POGS. 
Problems that took SDPT3 more than 150 seconds (of which there were many)
were aborted.

\subsection{Random problem classes}\label{sec:rand_probs}

We considered the following 9 problem classes: Basis pursuit, Entropy
maximization, Huber fitting, Lasso, Logistic regression, Linear programming,
Non-negative least-squares, Portfolio optimization, and Support vector machine
fitting. For each problem class, reasonable random instance were generated 
and solve; details about problem generation can be found in Appendix
\ref{sec:prob_gen}. For each problem class the number of non-zeros in $A$ was
varied on a logarithmic scale from 100 to 2 Billion. The aspect ratio of $A$
also varied from 1:10 to 10:1, with the orientation (wide or tall) chosen
depending on what was reasonable for each problem. We report running time
averaged over all aspect ratios.

The maximum number of iterations was set to $10^{4}$, but all problems
converged in fewer iterations, with most problems taking a couple of hundred
iterations. The relative tolerance was set to $10^{-3}$, and where solutions
from SDPT3 were available, we verified that the solutions produced by both
solvers matched to 3 decimal places. We omit SDPT3 running times for problems
involving exponential cones, since SDPT3 does not support them. 

\begin{figure}
\begin{center}
\includegraphics[scale=0.6]{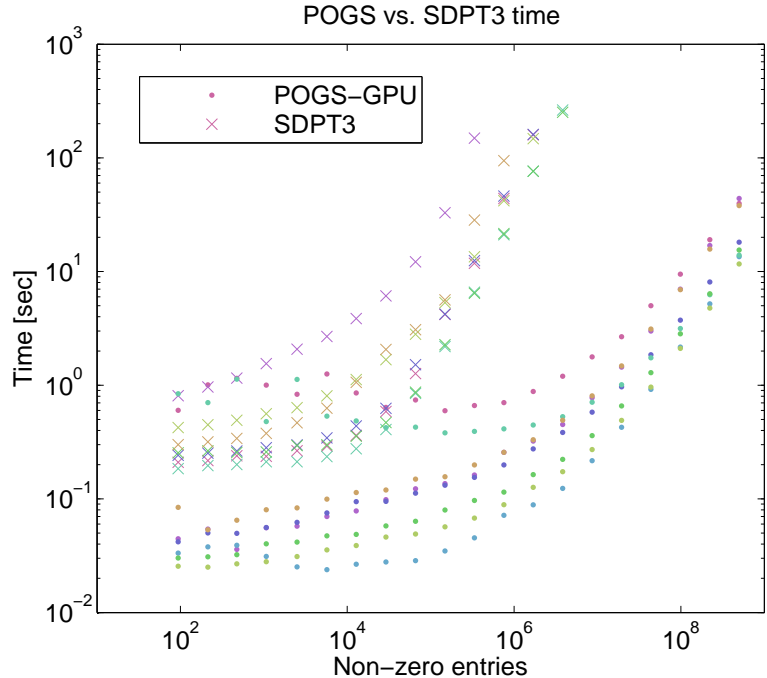}
\caption{POGS (GPU version) vs. SDPT3 for dense matrices (color represents problem class).}
\label{fig:pogs_dense}
\end{center}
\end{figure}

Figure \ref{fig:pogs_dense} compares the running time of POGS versus SDPT3, for problems where the constraint matrix $A$ is dense. 
We can make several general observations.
\begin{itemize}
\item POGS solves problems that are 3 orders of magnitude larger than SDPT3 in the same amount of time.
\item Problems that take 200 seconds in SDPT3 take 0.5 seconds in POGS.
\item POGS can solve problems with 2 Billion non-zeros in 10-50 seconds.
\item The variation in solve time across different problem classes
was similar for POGS and SDPT3, around one order of magnitude.
\end{itemize}
In summary, POGS is able to solve much larger problems, much faster
(to moderate precision).

\subsection{Radiation treatment planning}

Radiation treatment is used to radiate tumor cells in
cancer patients.
The goal of radiation treatment planning is to find a set of radiation beam
intensities that will deliver a specified radiation dosage to tumor cells,
while minimizing the impact on healthy cells. 
The problem can be stated directly in
graph form, with $x$ corresponding to the $n$ beam intensities to be found,
$y$ corresponding to the radiation dose received at the $m$ voxels,
and the matrix $A$ (whose elements are non-negative) giving the mapping from
the beams to the received dosages at the voxels.
This matrix comes from geometry, including 
radiation scattering inside the patient \cite{ahnesjo2006imrt}.
The objective $g$ is the indicator function of the non-negative orthant
(which imposes the constraint that $x_j \geq 0$), and $f$ is a separable function
of the form
\[
f_i(y_i) = \left\{ \begin{array}{ll}
w_i^+ y_i & \mbox{$i$ corresponds to a non-tumor voxel}\\
w_i^- \max (d_i-y_i,0)+
w_i^+ \max (y_i-d_i,0)
& \mbox{$i$ corresponds to a tumor voxel},
\end{array}\right.
\]
where $w_i^+ >0$ is the (given) weight associated with overdosing voxel $i$,
where $w_i^- >0$ is the (given) weight associated with underdosing voxel $i$,
and $d_i>0$ is the target dose, given for each tumor voxel.
We can also add the redundant constraint $y_i \geq 0$ by defining
$f_i(y_i)=\infty$ for $y_i<0$.

We present results for one instance of this problem, with $m=360000$ voxels 
and $n=360$ beams.  The matrix $A$ comes from a real patient, and the objective
parameters are chosen to achieve a good clinical plan.
The problem is small enough that it can be solved (to high accuracy) 
by an interior-point method, in around one hour.
POGS took a few seconds to solve the problem, producing a solution that was
extremely close to the one produced by the interior-point method.
In warm start mode, POGS could solve problem instances (obtained by varying the
objective parameters) in under one second, allowing for real-time
tuning of the treatment plan (by adjusting the objective function weights)
by a radiation oncologist.

\section{Acknowledgments}

We would like to thank Baris Ungun for testing POGS and providing valuable
feedback, as well as providing the radiation treatment data. We also thank
Michael Saunders for numerous discussions about solving large sparse systems.
This research was funded by DARPA XDATA and Adobe.

\appendix
\section{Problem generation details} \label{sec:prob_gen}
In this section we describe how the problems in \S\ref{sec:rand_probs} were generated. 

\subsection{Basis pursuit}
The basis pursuit problem \cite{chen1998atomic} seeks the smallest vector in the $\ell_1$-norm sense that satisfies a set of underdetermined linear equality 
constraints. The objective has the effect of finding a sparse solution. It can be stated as
\begin{align*}
	\begin{aligned}
	&\mini
	& & \|x\|_1 \\
	& \text{subject to} 
	& & b = A x,
	\end{aligned}
\end{align*}
with equivalent graph form representation
\begin{align*}
	\begin{aligned}
	&\mini
	& & I(y = b) + \|x\|_1 \\
	& \text{subject to} 
	& & y = A x.
	\end{aligned}
\end{align*}
The elements of $A$ were generated as $A_{ij} \sim \mathcal{N}(0, 1)$. To construct  $b$ we first generated a vector $v \in \mathbf{R}^n$  as
\[
	v_i \sim \left\{
	\begin{tabular}{ll}
		0 & with probability $p=1/2$ \\
		$\mathcal{N}(0, 1/n)$ & otherwise,
	\end{tabular}\right.
\]
we then let $b = Av$. In each instance we chose $m > n$.

\subsection{Entropy maximization}
The entropy maximization problem \cite{BoV:04} seeks a probability distribution with maximum entropy that satisfies a set of $m$ affine inequalities,
which can be interpreted as bounds on the expectations of arbitrary functions.
It can be stated as
\begin{align*}
	\begin{aligned}
	&\maxi
	& & -\textstyle{\sum}_{i=1}^nx_i\log x_i \\
	& \text{subject to} 
	& & \mathbf{1}^T x = 1, \quad A x \leq b,
	\end{aligned}
\end{align*}
with equivalent graph form representation
\begin{align*}
	\begin{aligned}
	&\mini
	& & I(y_{1:m} \leq b) + I(y_{m+1} = 1) + \textstyle{\sum}_{i=1}^nx_i\log x_i \\
	& \text{subject to} 
	& & y = \begin{bmatrix}A \\ \mathbf{1}^T\end{bmatrix} x.
	\end{aligned}
\end{align*}
The elements of $A$ were generated as $A_{ij} \sim \mathcal{N}(0,n)$. To construct $b$, we first generated a vector $v \in \mathbf{R}^n$ as $v_i \sim U[0, 1]$, then we set $b = Fv/(\mathbf{1}^Tv)$. This ensures that there exists a feasible $x$. In each instance we chose $m < n$.

\subsection{Huber fitting}
Huber fitting or robust regression \cite{huber1964robust} performs linear regression under the assumption that there are outliers in the data. The problem can be stated as
\begin{align*}
	\begin{aligned}
	&\mini
	& & \textstyle{\sum}_{i=1}^m\text{huber}(b_i - a_i^Tx),
	\end{aligned}
\end{align*}
where the Huber loss function is defined as
\[
	\text{huber}(x) = \left\{\begin{aligned} &(1/2)x^2 & |x| \leq 1 \\\ &|x| - (1/2) & |x| > 1 \end{aligned} \right.
\]
The graph form representation of this problem is
\begin{align*}
	\begin{aligned}
	&\mini &
	\textstyle{\sum}_{i=1}^n\text{huber}(b_i-y_i) \\
	& \text{subject to} 
	& y = Ax.
	\end{aligned}
\end{align*}
The elements of $A$ were generated as $A_{ij} \sim \mathcal{N}(0,n)$. To construct $b$, we first generated a vector $v \in \mathbf{R}^n$ as $v_i \sim \mathcal{N}(0,1/n)$
then we generated a noise vector $\varepsilon$ with elements
\[
	\varepsilon_i \sim \left\{
	\begin{tabular}{ll}
		$\mathcal{N}(0,1/4)$ & with probability $p=0.95$ \\
		$U[0, 10]$ & otherwise.
	\end{tabular}\right.
\]
Lastly we constructed $b = Av + \varepsilon$. In each instance we chose $m > n$.

\subsection{Lasso} \label{sec:apx_lasso}
The lasso problem \cite{tibshirani1996regression} seeks to perform linear regression under the assumption that the solution is sparse. An $\ell_1$ penalty is added to the objective to encourage sparsity. It can be stated as
\begin{align*}
	\begin{aligned}
	&\mini
	& & \|Ax-b\|_2^2 + \lambda \|x\|_1,
	\end{aligned}
\end{align*}
with graph form representation
\begin{align*}
	\begin{aligned}
	&\mini
	& & \|y-b\|_2 + \lambda\|x\|_1\\
	& \text{subject to} 
	& & y = Ax.
	\end{aligned}
\end{align*}
The elements of $A$ were generated as $A_{ij} \sim \mathcal{N}(0, 1)$. To construct  $b$ we first generated a vector $v \in \mathbf{R}^n$, with elements
\[
	v_i \sim \left\{
	\begin{tabular}{ll}
		0 & with probability $p=1/2$ \\
		$\mathcal{N}(0, 1/n)$ & otherwise.
	\end{tabular}\right.
\]
We then let $b = Av + \varepsilon$, where $\varepsilon$ represents the noise and was generated as $\varepsilon_i\sim\mathcal{N}(0,1/4)$. 
The value of $\lambda$ was set to $(1/5)\|A^Tb\|_\infty$.  This is a 
reasonable choice since
$\|A^Tb\|_\infty$ is the critical value of $\lambda$ above which the 
solution of the Lasso problem is $x=0$.
In each instance we chose $m < n$.

\subsection{Logistic regression}
Logistic regression \cite{hastie2009elements} fits a probability distribution to a binary class label. Similar to the Lasso problem (\ref{sec:apx_lasso}) a sparsifying $\ell_1$ penalty is often added to the coefficient vector. It can be stated as
\begin{align*}
	\begin{aligned}
	&\mini
	& &\textstyle{\sum}_{i=1}^m \left(\log(1+\exp (x^Ta_i)) - b_ix^Ta_i\right) + \lambda \|x\|_1,
	\end{aligned}
\end{align*}
where $b_i \in \{0,1\}$ is the class label of the $i$th sample, and $a_i^T$ is the 
$i$th row of $A$.
The graph form representation of this problem is
\begin{align*}
	\begin{aligned}
	&\mini
	& & \textstyle{\sum}_{i=1}^m \left(\log(1+\exp (y_i)) - b_iy_i\right) + \lambda \|x\|_1, \\
	& \text{subject to} 
	& & y = Ax.
	\end{aligned}
\end{align*}
The elements of $A$ were generated as $A_{ij} \sim \mathcal{N}(0, 1)$. To construct $b$ we first generated a vector $v \in \mathbf{R}^n$, with elements
\[
	v_i \sim \left\{
	\begin{tabular}{ll}
		0 & with probability $p=1/2$ \\
		$\mathcal{N}(0, 1/n)$ & otherwise.
	\end{tabular}\right.
\]
We then constructed the entries of $b$ as
\[
	b_i \sim \left\{
	\begin{tabular}{ll}
		0 & with probability $p=1/(1+ \exp(-a_i^Tv))$ \\
		1 & otherwise.
	\end{tabular}\right.
\]
The value of $\lambda$ was set to $(1/10)\|A^T((1/2)\mathbf{1} - b)\|_\infty$.
($\|A^T((1/2)\mathbf{1} - b)\|_\infty$ is the critical of 
$\lambda$ above which the solution is $x=0$.)
  In each instance we chose $m > n$.

\subsection{Linear program}
Linear programs \cite{BoV:04} seek to minimize a linear function subject to linear inequality constraints. It can be stated as
\[
	\begin{aligned}
	&\mini
	& & c^Tx \\
	& \text{subject to} 
	& & Ax \leq b,
	\end{aligned}
\]
and has graph form representation
\[
	\begin{aligned}
	&\mini
	& & c^Tx + I(y \leq b) \\
	& \text{subject to} 
	& & y = Ax.
	\end{aligned}
\]
The elements of $A$ were generated as $A_{ij} \sim \mathcal{N}(0, 1)$. To construct $b$ we first generated a vector $v \in \mathbf{R}^n$, with elements
\[
	v_i \sim \mathcal{N}(0, 1/n).
\]
We then generated $b$ as $b = Av + \varepsilon$, where $\varepsilon_i \sim U[0, 1/10]$. The vector $c$ was constructed in a similar fashion. First we generate a vector $u \in \mathbf{R}^m$, with elements
\[
	u_i \sim U[0,1],
\]
then we constructed $c = -A^Tu$. This method guarantees that the problem is bounded. In each instance we chose $m > n$.

\subsection{Non-negative least-squares}
Non-negative least-squares \cite{chen2009nonnegativity} seeks a minimizer of 
a least-squares problem subject to the solution vector being non-negative. This comes up in applications where the solution represents real quantities. The problem can be stated as
\[
	\begin{aligned}
	&\mini
	& & \|Ax-b\|_2 \\
	& \text{subject to} 
	& & x \geq 0,
	\end{aligned}
\]
and has graph form representation
\[
	\begin{aligned}
	&\mini
	& & \|y-b\|_2^2 + I(x\geq 0) \\
	& \text{subject to} 
	& & y = Ax.
	\end{aligned}
\]
The elements of $A$ were generated as $A_{ij} \sim \mathcal{N}(0, 1)$. To construct $b$ we first generated a vector $v \in \mathbf{R}^n$, with elements
\[
	v_i \sim \mathcal{N}(1/n, 1/n).
\]
We then generated $b$ as $b = Av + \varepsilon$, where $\varepsilon_i \sim \mathcal{N}(0,1/4)$. In each instance we chose $m > n$.

\subsection{Portfolio optimization}
Portfolio optimization or optimal asset allocation seeks to maximize the risk adjusted return of a portfolio. A common assumption is the $k$-factor risk model \cite{connor1993arbitrage}, which states that the return covariance matrix is the sum of a diagonal plus a rank $k$ matrix. The problem can be stated as
\[
	\begin{aligned}
	&\maxi
	& & \mu^Tx- \gamma x^T(FF^T + D)x \\
	& \text{subject to} 
	& & x \geq 0, \quad \mathbf{1}^Tx = 1
	\end{aligned}
\]
where $F \in \reals^{n \times k}$ and $D$ is diagonal. An equivalent graph form representation is given by
\[
	\begin{aligned}
	&\mini
	& &-x^T\mu + \gamma x^TDx + I(x \geq 0)+ \gamma y_{1:m}^Ty_{1:m} + I(y_{m+1} = 1)\\
	& \text{subject to} 
	& & y = \begin{bmatrix}
		F^T \\ \mathbf{1}^T
	\end{bmatrix}x.
	\end{aligned}
\]
The elements of $A$ were generated as $A_{ij} \sim \mathcal{N}(0,1)$. The diagonal of $D$ was generated as $D_{ii} \sim U[0,\sqrt{k}]$ and the the mean return $\mu$ was generated as $\mu_i \sim \mathcal{N}(0, 1)$. The risk aversion factor $\gamma$ was set to 1. In each instance we chose $n > k$.

\subsection{Support vector machine}
The support vector machine \cite{cortes1995support} problem seeks a separating hyperplane classifier for a problem with two classes. The problem can be stated as
\[
	\begin{aligned}
	&\mini
	& & x^Tx + \lambda \textstyle{\sum}_{i = 1}^m\max(0, b_ia_i^Tx + 1), 
	\end{aligned}
\]
where $b_i \in \{-1, +1\}$ is a class label and $a_i^T$ is the $i$th row of $A$.
It has graph form representation
\[
	\begin{aligned}
	&\mini
	& & \lambda \textstyle{\sum}_{i = 1}^m\max(0, y_i + 1) + x^Tx\\
	& \text{subject to} 
	& & y = \mathbf{diag}(b)Ax.
	\end{aligned}
\]
The vector $b$ was chosen to so that the first $m/2$ elements belong to one class and the second $m/2$ belong to the other class. Specifically 
\[
	b_i = \left\{\begin{tabular}{ll}
		$+1$ & $i \leq m/2$ \\
		$-1$ & otherwise.
	\end{tabular}\right.
\]
Similarly, the elements of $A$ were generated as
\[
	A_{ij} \sim \left\{\begin{tabular}{ll}
		$\mathcal{N}(+1/n, 1/n)$ & $i \leq m/2$ \\
		$\mathcal{N}(-1/n, 1/n)$ & otherwise.
	\end{tabular}\right.
\]
This choice of $A$ causes the rows of $A$ to form two distinct clusters. In each instance we chose $m > n$.

\newpage
\bibliography{pogs}

\newcommand{\etalchar}[1]{$^{#1}$}
\begin{thebibliography}{BEGFB94}

\bibitem[AHIM06]{ahnesjo2006imrt}
A.~Ahnesj{\"o}, B.~H{\aa}rdemark, U.~Isacsson, and A.~Montelius.
\newblock The {IMRT} information process --- {M}astering the degrees of freedom
  in external beam therapy.
\newblock {\em Physics in Medicine and Biology}, 51(13):R381--R402, 2006.

\bibitem[AHW12]{annergren2012admm}
M.~Annergren, A.~Hansson, and B.~Wahlberg.
\newblock An {ADMM} algorithm for solving $\ell_1$ regularized {MPC}.
\newblock {\em arXiv preprint arXiv:1203.4070}, 2012.

\bibitem[BEGFB94]{boyd1994linear}
S.~Boyd, L.~El~Ghaoui, E.~Feron, and V.~Balakrishnan.
\newblock {\em Linear matrix inequalities in system and control theory},
  volume~15.
\newblock SIAM, 1994.

\bibitem[BMOW13]{boyd2013performance}
S.~Boyd, M.~Mueller, B.~O'Donoghue, and Y.~Wang.
\newblock Performance bounds and suboptimal policies for multi-period
  investment.
\newblock {\em Foundations and Trends in Optimization}, 1(1):1--69, 2013.

\bibitem[BPC{\etalchar{+}}11]{boyd2011distributed}
S.~Boyd, N.~Parikh, E.~Chu, B.~Peleato, and J.~Eckstein.
\newblock Distributed optimization and statistical learning via the alternating
  direction method of multipliers.
\newblock {\em Foundations and Trends in Machine Learning}, 3(1):1--122, 2011.

\bibitem[Bra10]{bradley2010algorithms}
A.~M. Bradley.
\newblock {\em Algorithms for the equilibration of matrices and their
  application to limited-memory quasi-{N}ewton methods}.
\newblock PhD thesis, Stanford University, 2010.

\bibitem[BTN01]{ben2001lectures}
A.~Ben-Tal and A.~Nemirovski.
\newblock {\em Lectures on modern convex optimization: {A}nalysis, algorithms,
  and engineering applications}, volume~2.
\newblock SIAM, 2001.

\bibitem[BV04]{BoV:04}
S.~Boyd and L.~Vandenberghe.
\newblock {\em Convex Optimization}.
\newblock Cambridge University Press, 2004.

\bibitem[CDHR08]{chen2008algorithm}
Y.~Chen, T.~A. Davis, W.~W. Hager, and S.~Rajamanickam.
\newblock Algorithm 887: {CHOLMOD}, supernodal sparse {C}holesky factorization
  and update/downdate.
\newblock {\em ACM Transactions on Mathematical Software}, 35(3):22, 2008.

\bibitem[CDS98]{chen1998atomic}
S.~S. Chen, D.~L. Donoho, and M.~A. Saunders.
\newblock Atomic decomposition by basis pursuit.
\newblock {\em SIAM Journal on Scientific Computing}, 20(1):33--61, 1998.

\bibitem[CHW{\etalchar{+}}13]{coates2013deep}
A.~Coates, B.~Huval, T.~Wang, D.~Wu, B.~Catanzaro, and A.~Y. Ng.
\newblock Deep learning with {COTS HPC} systems.
\newblock In {\em Proceedings of the 30th International Conference on Machine
  Learning}, pages 1337--1345, 2013.

\bibitem[CK93]{connor1993arbitrage}
G.~Connor and R.~A. Korajczyk.
\newblock The arbitrage pricing theory and multifactor models of asset returns.
\newblock {\em Handbooks in Operations Research and Management Science}, 9,
  1993.

\bibitem[CM87]{calamai1987projected}
P.~H. Calamai and J.~J. Mor{\'e}.
\newblock Projected gradient methods for linearly constrained problems.
\newblock {\em Mathematical Programming}, 39(1):93--116, 1987.

\bibitem[CMM{\etalchar{+}}11]{ciresan2011flexible}
D.~C. Ciresan, U.~Meier, J.~Masci, L.~M. Gambardella, and J.~Schmidhuber.
\newblock Flexible, high performance convolutional neural networks for image
  classification.
\newblock {\em International Joint Conference on Artificial Intelligence},
  22(1):1237--1242, 2011.

\bibitem[COPB13]{chu2013primal}
E.~Chu, B.~O'Donoghue, N.~Parikh, and S.~Boyd.
\newblock A primal-dual operator splitting method for conic optimization.
\newblock 2013.

\bibitem[CP09]{chen2009nonnegativity}
D.~Chen and R.~J. Plemmons.
\newblock Nonnegativity constraints in numerical analysis.
\newblock pages 109--140, 2009.

\bibitem[CP11]{chambolle2011first}
A.~Chambolle and T.~Pock.
\newblock A first-order primal-dual algorithm for convex problems with
  applications to imaging.
\newblock {\em Journal of Mathematical Imaging and Vision}, 40(1):120--145,
  2011.

\bibitem[CV95]{cortes1995support}
C.~Cortes and V.~Vapnik.
\newblock Support-vector networks.
\newblock {\em Machine Learning}, 20(3):273--297, 1995.

\bibitem[DR56]{douglas1956numerical}
J.~Douglas and H.~H. Rachford.
\newblock On the numerical solution of heat conduction problems in two and
  three space variables.
\newblock {\em Transactions of the American Mathematical Society}, pages
  421--439, 1956.

\bibitem[DY14]{davis2014convergence}
D.~Davis and W.~Yin.
\newblock Convergence rate analysis of several splitting schemes.
\newblock {\em arXiv preprint arXiv:1406.4834}, 2014.

\bibitem[EB92]{eckstein1992douglas}
J.~Eckstein and D.~P. Bertsekas.
\newblock On the {D}ouglas-{R}achford splitting method and the proximal point
  algorithm for maximal monotone operators.
\newblock {\em Mathematical Programming}, 55(1-3):293--318, 1992.

\bibitem[ES08]{eckstein2008family}
J.~Eckstein and B.~F. Svaiter.
\newblock A family of projective splitting methods for the sum of two maximal
  monotone operators.
\newblock {\em Mathematical Programming}, 111(1-2):173--199, 2008.

\bibitem[GB14a]{giselsson2014diagonal}
P.~Giselsson and S.~Boyd.
\newblock Diagonal scaling in {D}ouglas-{R}achford splitting and {ADMM}.
\newblock In {\em 53rd IEEE Conference on Decision and Control}, 2014.

\bibitem[GB14b]{giselsson2014metric}
P.~Giselsson and S.~Boyd.
\newblock Metric selection in {D}ouglas-{R}achford splitting and {ADMM}.
\newblock {\em arXiv preprint arXiv:1410.8479}, 2014.

\bibitem[GB14c]{giselsson2014preconditioning}
P.~Giselsson and S.~Boyd.
\newblock Preconditioning in fast dual gradient methods.
\newblock {\em 53rd IEEE Conference on Decision and Control}, 2014.

\bibitem[Gis15]{giselsson2015tight}
P.~Giselsson.
\newblock Tight linear convergence rate bounds for {D}ouglas-{R}achford
  splitting and {ADMM}.
\newblock {\em arXiv preprint arXiv:1503.00887}, 2015.

\bibitem[GM75]{Glowinski1975}
R.~Glowinski and A.~Marroco.
\newblock Sur l'approximation, par \'el\'ments finis d'ordre un, et la
  r\'esolution, par p\'enalisation-dualit\'e d'une classe de probl\'emes de
  {D}irichlet non lin\'eaires.
\newblock {\em Mathematical Modelling and Numerical Analysis}, 9(R2):41--76,
  1975.

\bibitem[GOSB14]{goldstein2014fast}
T.~Goldstein, B.~O'Donoghue, S.~Setzer, and R.~Baraniuk.
\newblock Fast alternating direction optimization methods.
\newblock {\em SIAM Journal on Imaging Sciences}, 7(3):1588--1623, 2014.

\bibitem[GTSJ13]{ghadimi2013optimal}
E.~Ghadimi, A.~Teixeira, I.~Shames, and M.~Johansson.
\newblock Optimal parameter selection for the alternating direction method of
  multipliers ({ADMM}): quadratic problems.
\newblock {\em IEEE Transactions on Automatic Control}, 60:644--658, 2013.

\bibitem[HS52]{hestenes1952methods}
M.~R. Hestenes and E.~Stiefel.
\newblock Methods of conjugate gradients for solving linear systems.
\newblock {\em Joural of Research of the National Bureau of Standards},
  49(6):409--436, 1952.

\bibitem[HTF09]{hastie2009elements}
T.~Hastie, R.~Tibshirani, and T.~Friedman.
\newblock {\em The elements of statistical learning}, volume~2.
\newblock Springer, 2009.

\bibitem[Hub64]{huber1964robust}
P.~J. Huber.
\newblock Robust estimation of a location parameter.
\newblock {\em The Annals of Mathematical Statistics}, 35(1):73--101, 1964.

\bibitem[HYW00]{he2000alternating}
B.~S. He, H.~Yang, and S.~L. Wang.
\newblock Alternating direction method with self-adaptive penalty parameters
  for monotone variational inequalities.
\newblock {\em Journal of Optimization Theory and Applications},
  106(2):337--356, 2000.

\bibitem[KSH12]{krizhevsky2012imagenet}
A.~Krizhevsky, I.~Sutskever, and G.~E. Hinton.
\newblock Imagenet classification with deep convolutional neural networks.
\newblock In {\em Advances in Neural Information Processing Systems}, pages
  1097--1105, 2012.

\bibitem[LM79]{lions1979splitting}
P.~L. Lions and B.~Mercier.
\newblock Splitting algorithms for the sum of two nonlinear operators.
\newblock {\em SIAM Journal on Numerical Analysis}, 16(6):964--979, 1979.

\bibitem[NCL{\etalchar{+}}11]{ngiam2011optimization}
J.~Ngiam, A.~Coates, A.~Lahiri, B.~Prochnow, Q.~V. Le, and A.~Y. Ng.
\newblock On optimization methods for deep learning.
\newblock In {\em Proceedings of the 28th International Conference on Machine
  Learning}, pages 265--272, 2011.

\bibitem[NLR{\etalchar{+}}15]{nishihara2015general}
R.~Nishihara, L.~Lessard, B.~Recht, A.~Packard, and M.~I. Jordan.
\newblock A general analysis of the convergence of {ADMM}.
\newblock {\em arXiv preprint arXiv:1502.02009}, 2015.

\bibitem[NW99]{wrightnocedal1999numerical}
J.~Nocedal and S.~Wright.
\newblock {\em Numerical Optimization}, volume~2.
\newblock Springer, 1999.

\bibitem[OSB13]{o2013splitting}
B.~O'Donoghue, G.~Stathopoulos, and S.~Boyd.
\newblock A splitting method for optimal control.
\newblock {\em IEEE Transactions on Control Systems Technology},
  21(6):2432--2442, 2013.

\bibitem[OV14]{o2014primal}
D.~O'Connor and L.~Vandenberghe.
\newblock Primal-dual decomposition by operator splitting and applications to
  image deblurring.
\newblock {\em SIAM Journal on Imaging Sciences}, 7(3):1724--1754, 2014.

\bibitem[OW06]{olafsson2006efficient}
A.~Olafsson and S.~Wright.
\newblock Efficient schemes for robust {IMRT} treatment planning.
\newblock {\em Physics in Medicine and Biology}, 51(21):5621--5642, 2006.

\bibitem[PB13a]{parikh2013block}
N.~Parikh and S.~Boyd.
\newblock Block splitting for distributed optimization.
\newblock {\em Mathematical Programming Computation}, pages 1--26, 2013.

\bibitem[PB13b]{parikh2013proximal}
N.~Parikh and S.~Boyd.
\newblock Proximal algorithms.
\newblock {\em Foundations and Trends in Optimization}, 1(3):123--231, 2013.

\bibitem[PC11]{pock2011diagonal}
T.~Pock and A.~Chambolle.
\newblock Diagonal preconditioning for first order primal-dual algorithms in
  convex optimization.
\newblock In {\em IEEE International Conference on Computer Vision}, pages
  1762--1769, 2011.

\bibitem[Pol87]{polyak1987introduction}
B.~Polyak.
\newblock Introduction to optimization.
\newblock {\em Optimization Software Inc., Publications Division, New York},
  1987.

\bibitem[PS75]{paige1975solution}
C.~C. Paige and M.~A. Saunders.
\newblock Solution of sparse indefinite systems of linear equations.
\newblock {\em SIAM Journal on Numerical Analysis}, 12(4):617--629, 1975.

\bibitem[PS82]{paige1982lsqr}
C.~C. Paige and M.~A. Saunders.
\newblock {LSQR}: {A}n algorithm for sparse linear equations and sparse least
  squares.
\newblock {\em ACM Transactions on Mathematical Software}, 8(1):43--71, 1982.

\bibitem[Rui01]{ruiz2001scaling}
D.~Ruiz.
\newblock A scaling algorithm to equilibrate both rows and columns norms in
  matrices.
\newblock Technical report, Rutherford Appleton Laboratory, 2001.
\newblock Technical Report RAL-TR-2001-034.

\bibitem[Sho98]{shor1998nondifferentiable}
N.~Z. Shor.
\newblock {\em Nondifferentiable optimization and polynomial problems}.
\newblock Kluwer Academic Publishers, 1998.

\bibitem[SK67]{sinkhorn1967concerning}
R.~Sinkhorn and P.~Knopp.
\newblock Concerning nonnegative matrices and doubly stochastic matrices.
\newblock {\em Pacific Journal of Mathematics}, 21(2):343--348, 1967.

\bibitem[Spi85]{spingarn1985applications}
J.~E. Spingarn.
\newblock Applications of the method of partial inverses to convex programming:
  decomposition.
\newblock {\em Mathematical Programming}, 32(2):199--223, 1985.

\bibitem[Tib96]{tibshirani1996regression}
R.~Tibshirani.
\newblock Regression shrinkage and selection via the lasso.
\newblock {\em Journal of the Royal Statistical Society}, pages 267--288, 1996.

\bibitem[TTT99]{toh1999sdpt3}
K.~Toh, M.~J. Todd, and R.~H. T{\"u}t{\"u}nc{\"u}.
\newblock {SDPT3}-a {MATLAB} software package for semidefinite programming,
  version 1.3.
\newblock {\em Optimization Methods and Software}, 11(1-4):545--581, 1999.

\bibitem[Van95]{vanderbei1995symmetric}
R.~J. Vanderbei.
\newblock Symmetric quasidefinite matrices.
\newblock {\em SIAM Journal on Optimization}, 5(1):100--113, 1995.

\bibitem[WB14]{wang2014bregman}
H.~Wang and A.~Banerjee.
\newblock Bregman alternating direction method of multipliers.
\newblock In {\em Advances in Neural Information Processing Systems}, pages
  2816--2824, 2014.

\end{thebibliography}

\end{document}